\documentclass[a4paper,11pt]{article}
\usepackage{a4wide}
\usepackage{theorem}
\usepackage{amsmath}
\usepackage{array}
\usepackage{amssymb}
\usepackage{amsfonts}
\usepackage[french,english]{babel}
\usepackage{epsf}
\usepackage{epsfig}

\newtheorem{theo}{\indent Theorem\newline}[section]

{\theorembodyfont{\rmfamily}%
\newtheorem{rem}[theo]{\noindent Remark}}
{\theorembodyfont{\rmfamily}%
 \theoremstyle{break}%
}
\newtheorem{prop}[theo]{\indent Proposition\newline}
\newtheorem{lemma}[theo]{\indent Lemma\newline}
\newtheorem{cor}[theo]{\indent Corollary\newline}

 \def\N{{\Bbb{N}}}
\def\Z{{\Bbb{Z}}}

\def\R{{\Bbb{R}}}
\def\C{{\Bbb{C}}}

\usepackage{eufrak}

\newcommand{\mult}{\mathop{\rm mult}\nolimits}
\newcommand{\ind}{\mathop{\rm ind}\nolimits}

\newcommand{\coker}{\mathop{\rm coker}\nolimits}

\newlength{\indentation}%
\makeatletter
\newcommand\@makefntextsans[1]{%
    \parindent 0em%
    \noindent%
    \hb@xt@0em{\hss}%
    #1}
\def\footnotetextsans{%
     \@ifnextchar [\@xfootnotenextsans%
       {\@footnotetextsans}}
\def\@xfootnotenextsans[#1]{%
  \begingroup%
     \csname c@\@mpfn\endcsname #1\relax%
  \endgroup%
  \@footnotetextsans}
\long\def\@footnotetextsans#1{\insert\footins{%
    \reset@font\footnotesize%
    \interlinepenalty\interfootnotelinepenalty%
    \splittopskip\footnotesep%
    \splitmaxdepth \dp\strutbox \floatingpenalty \@MM%
    \hsize\columnwidth \@parboxrestore%
    \color@begingroup%
      \@makefntextsans{%
        \rule\z@\footnotesep\ignorespaces#1\@finalstrut\strutbox}
    \color@endgroup}}
\makeatother

\begin{document}

\cleardoublepage
\title{Towards relative invariants of real symplectic $4$-manifolds}
\author{Jean-Yves Welschinger}
\maketitle

\makeatletter\renewcommand{\@makefnmark}{}\makeatother
\footnotetextsans{Keywords : Real symplectic manifold, 
rational curve, enumerative geometry.}
\footnotetextsans{AMS Classification : 53D45, 14N35, 14N10, 14P99.}

{\bf Abstract:}

Let $(X , \omega , c_X)$ be a real symplectic $4$-manifold with real part 
$\R X$. Let $L \subset \R X$ be a smooth curve such that $[L] = 0 \in
H_1 (\R X ; \Z / 2\Z)$. We construct invariants under deformation of
the quadruple $(X , \omega , c_X , L)$ by counting the number of real
rational $J$-holomorphic curves which realize a given homology class $d$,
pass through an appropriate number of points and are tangent to $L$. As 
an application, we prove a relation between the count of real rational 
$J$-holomorphic curves done in \cite{Wels1} and the count of
reducible real rational curves done in \cite{Wels2}. Finally, we show how 
these techniques also allow to extract an integer valued invariant from a classical 
problem of real enumerative geometry, namely about counting the number of real plane conics
tangent to five given generic real conics.

\section{Statement of the results}

Let $(X , \omega , c_X)$ be a {\it real symplectic $4$-manifold}, that is a triple made of
a smooth compact $4$-manifold $X$, a symplectic form $\omega$ on $X$ and an involution $c_X$ on $X$
such that $c_X^* \omega =- \omega$. The fixed 
point set of $c_X$ is called {\it the real part of $X$} and is
denoted by $\R X$. It is assumed to be non empty here so that it is a smooth lagrangian 
surface of $(X, \omega)$. We label its connected components by $(\R X)_1 , \dots ,
(\R X)_N$. Let $L \subset \R X$ be a smooth curve which represents $0$ in $H_1 (\R X ; \Z / 2\Z)$,
and $B \subset \R X$ be a surface having $L$ as a boundary.

\subsection{Definitions}

Let $l \gg 1$ be an integer large enough and ${\cal J}_\omega$ be the space of 
almost complex structures of $X$ which are tamed by $\omega$ and of class $C^l$.
Let $\R {\cal J}_\omega$ be the subspace of ${\cal J}_\omega$ made of almost complex 
structures $J$ for which the involution $c_X$ is $J$-antiholomorphic. These two spaces
are separable Banach manifolds which are non empty and contractible (see \S $1.1$
of \cite{Wels1} for the real case). Assume that the first Chern class $c_1 (X)$ of the symplectic 
$4$-manifold $(X , \omega)$ is not a torsion element
in $H^2 (X ; \Z)$ and let $d \in H_2 (X ; \Z)$ be
a homology class satisfying $c_1 (X) d > 1$, $c_1 (X) d \neq 4$ and $(c_X)_* d = -d$.
Let $\underline{x} = (x_1 , \dots , x_{c_1 (X) d -2}) \in X^{c_1 (X) d -2}$ be a 
{\it real configuration} of $c_1 (X) d -2$ distinct points of $X$, that is an ordered 
subset of distinct points of $X$ which is globally invariant under $c_X$. For
$j \in \{ 1 , \dots , N \}$, we set $r_j = \# ( \underline{x} \cap (\R X)_j)$ and
$r = (r_1 , \dots , r_N)$ so that the $N$-tuple $r$ encodes the equivariant isotopy
class of $\underline{x}$. We will assume throughout the paper that $r \neq (0 , \dots , 0)$.
Finally, denote by $I$ the subset of those $i \in \{ 1 , \dots,
c_1 (X) d -2 \}$ for which $x_i$ is fixed by the involution $c_X$, so that 
$I \neq \emptyset$. For each 
$i \in I$, choose a line $T_i$ in the tangent plane $T_{x_i} \R X$.
Let $J \in \R {\cal J}_\omega$ be generic enough. Then, as in \cite{Wels2}, we denote by
${\cal C}usp^d (J , \underline{x})$ (resp. ${\cal R}ed^d (J , \underline{x})$,
${\cal T}an^d (J , \underline{x})$) the finite set of real rational cuspidal (resp.
reducible, whose tangent line at some point $x_i$, $i \in I$, is $T_i$)
$J$-holomorphic curves which realize the homology class $d$ and pass
through $\underline{x}$. Likewise, we denote by ${\cal T}an^d_L (J , \underline{x})$
the finite set of real rational $J$-holomorphic curves which realize the homology class $d$, pass
through $\underline{x}$ and are tangent to $L$. Note that the genericity assumption
on $J \in \R {\cal J}_\omega$ implies that the non-trivial point of contact of the curve 
with $L$ is unique and of order two.
Also, all these curves have only transversal double points as singularities lying outside 
of $\underline{x}$,
with the exception of elements of ${\cal C}usp^d (J , \underline{x})$ which have 
in addition a unique real ordinary cusp. Let
$C \in {\cal T}an^d_L (J , \underline{x}) \cup {\cal C}usp^d (J , \underline{x}) \cup 
{\cal R}ed^d (J , \underline{x}) \cup
{\cal T}an^d (J , \underline{x})$, we define the {\it mass} of 
$C$ and denote by $m(C)$ its number of real isolated double points. Here, a real double
point is said to be {\it isolated} when it is the local intersection of two complex 
conjugated branches, whereas it is said to be {\it non isolated} when it is
the local intersection of two real branches. Let $C \in {\cal T}an^d_L (J , \underline{x})$
and $y$ be its point of contact with $L$. Then, either $\R C$ is locally
included in $B$ near $y$, or its intersection with $B$ is locally restricted to $\{ y \}$. 
We define the {\it contact
index} $<C,B>$ to be $-1$ in the first case and $+1$ in the second. Likewise, if
$C \in {\cal C}usp^d (J , \underline{x})$ (resp. $C \in {\cal T}an^d (J , \underline{x})$), 
then its cuspidal point (resp. its tangent line $T_i$, $i \in I$) is unique and we define
$<C,B>$ to be $-1$ if it is outside $B$ or $+1$ if it is inside. Finally,
if $C$ belongs to ${\cal R}ed^d (J , \underline{x})$ and $C_1$, $C_2$
denote its irreducible components, then both these components are reals and we set
$$\mult_B (C) = \sum_{y \in \R C_1 \cap \R C_2} <y,B>,$$
where $<y,B>$ equals $-1$ if $y$ is outside $B$ or $+1$ if it is inside.

\subsection{Statement of the results}
\label{subsectresults}

We set:
$$\Gamma^{d,B}_r (J , \underline{x}) = \sum_{C \in 
\cup {\cal T}an^d_L (J , \underline{x}) \cup {\cal T}an^d (J , \underline{x})
\cup {\cal C}usp^d (J , \underline{x})
}
(-1)^{m(C)} <C,B> - \sum_{C \in {\cal R}ed^d (J , \underline{x})} (-1)^{m(C)} \mult_B (C).$$
\begin{theo}
\label{maintheo}
Let $(X, \omega , c_X)$ be a real symplectic $4$-manifold and $B \subset \R X$ be
a surface with boundary $L$. The connected components of $\R X$ are labelled by 
$(\R X)_1 , \dots , (\R X)_N$. Let $d \in H_2 (X ; \Z)$ be
such that $c_1 (X) d > 1$ and $c_1 (X) d \neq 4$, and $\underline{x} \subset X \setminus L$ be 
a real configuration of
$c_1 (X) d - 2$ distinct points. For $j \in \{ 1 , \dots , N \}$, denote by $r_j$ 
the cardinality of $\underline{x} \cap
(\R X)_j$ and by $r = (r_1 , \dots , r_N)$, which is assumed to be different from
$(0 , \dots , 0)$. Finally, let $J \in \R {\cal J}_\omega$
be generic enough so that the integer $\Gamma^{d,B}_r (J , \underline{x})$ is well defined.
Then, this integer $\Gamma^{d,B}_r (J , \underline{x})$ neither depends on the choice of
$J$, nor on the choice of $\underline{x}$.
\end{theo}
From this theorem, the integer $\Gamma^{d,B}_r (J , \underline{x})$ can 
be denoted without ambiguity by $\Gamma^{d,B}_r$, and when it is not well defined,
we set $\Gamma^{d,B}_r = 0$. Note that the condition $c_1 (X) d \neq 4$ is to avoid appearance 
of multiple curves, see Remark $1.11$ of \cite{Wels2}.
\begin{rem}
1) In particular, the integer $\Gamma^{d,B}_r (J , \underline{x})$ does not depend on the
relative position of $\underline{x}$ with respect to $B$, it only depends on $r$.

2) When $B = \emptyset$, $\Gamma^{d,B}_r = - \Gamma^d_r$, where $\Gamma^d_r$ is the invariant
defined in \cite{Wels2}. Theorem \ref{maintheo} then follows
from Theorem $0.1$ of \cite{Wels2}.

3) One has $\Gamma^{d,B}_r = - \Gamma^{d,\R X \setminus B}_r$.
\end{rem}

We denote by
$\Gamma^{d,B} (T)$ the generating function $\sum_{r \in \N^N} \Gamma^{d,B}_r T^r \in \Z [T_1,
\dots , T_N]$, where $T^r =T_1^{r_1} \dots T_N^{r_N}$. This polynomial function is of
the same parity as $c_1 (X)d$ and each of its monomial actually only depends on one
indeterminate. It follows from Theorem \ref{maintheo} that the
function $\Gamma^B : d \in H_2 (X ; \Z) \mapsto \Gamma^{d,B} (T) \in \Z [T]$ only depends on
the quadruple $(X , \omega , c_X , B)$. Moreover, it is invariant under
deformation of this quadruple, that is if $\omega_t$ is a continuous
family of symplectic forms on $X$ for which $c_X^* \omega_t = - \omega_t$ and $B_t$ is an isotopy
of compact surfaces in $\R X$, then
this function is the same for all $(X ,\omega_t , c_X , B_t)$. 
\begin{theo}
\label{theorelation}
Under the hypothesis of Theorem \ref{maintheo}, assume that $B$ is a disk in $\R X$. Then
$2\chi^d_{r+1} = \Gamma^{d,B}_r + \Gamma^d_r$. Moreover,

1) If $(X, \omega , c_X)$ is the complex projective plane equipped with its standard
symplectic form and real structure, then $\Gamma^{d,B}_r = \Gamma^d_r$. 

2) If $(X, \omega , c_X)$ is the hyperboloid $(\C P^1 \times \C P^1, \omega_{\C P^1} \oplus 
\omega_{\C P^1}, conj \times conj)$, then $\Gamma^{d,B}_r = 2\chi^d_{r+1} + \Gamma^d_r$.
\end{theo}
(Remember that the integer $\chi^d_{r+1}$ has been defined in \cite{Wels0}, \cite{Wels1}
and the integer $\Gamma^d_r$ in \cite{Wels2}. Note that when $\R X$ is connected, $r \in \N^*$.)
\begin{cor}
Under the hypothesis of Theorem \ref{theorelation}, we have
$\chi^d_{r+1} = \Gamma^d_r =  \Gamma^{d,B}_r$ in the case of the complex projective plane
and $\Gamma^{d,B}_r = 2 \chi^d_{r+1}$, $\Gamma^d_r = 0$ in the case of the hyperboloid. $\square$
\end{cor}
The first equality of this corollary has been announced in \cite{Wels2}, Proposition $0.3$. 
It provides a relation between the count of real rational $J$-holomorphic curves done in 
\cite{Wels1} and the count of reducible and cuspidal curves done in \cite{Wels2}. Does such
a relation have a complex analog? 

\subsection{More tangency conditions, the case of conics}
\label{subsectcon}

It is possible to extend the above results to curves having more than one tangency condition with $L$,
at least in the case of plane conics (see also \S \ref{subsectfinal}). We illustrate this phenomenum 
here on the following classical
problem of real enumerative geometry, solved by De Joncqui\`eres in $1859$: there are $3264$ conics
which are tangent to five given generic conics in the complex projective plane. If the five given 
conics are
real, then the number of real conics tangent to them of course depends on the choice of the conics.
We however show here how it is possible to extract an integer valued invariant from this problem.

Let $B_1 , \dots , B_5$ be five embedded disks in $\R P^2$ which are transversal to each other and
$L_i = \partial B_i$, $i \in \{ 1, \dots , 5 \}$. Let $J \in \R {\cal J}_\omega$ be generic enough.
We denote by ${\cal C}on (J)$ the finite set of real conics tangent to $L_1 , \dots , L_5$ and
by ${\cal C}on_{red} (J)$ the finite set of real reducible conics, that is pairs of 
$J$-holomorphic lines, tangent 
to four out of these five curves $L_1 , \dots , L_5$. Let $C \in {\cal C}on (J)$, we set
$<C,B> = \Pi_{i=1}^5 < C, B_i >$. In the same way, let $C \in {\cal C}on_{red} (J)$ and $i_1 , \dots ,
i_4 \in \{ 1, \dots , 5 \}$ be such that $C$ is tangent to $L_{i_1} , \dots , L_{i_4}$. We set
$<C ,  B> = \Pi_{j=1}^4 < C, B_{i_j} >$ and $\mult_B (C) = +1$ if the singular point of $C$ belongs
to $B_{i_5}$ and $-1$ otherwise. Set
$$\Gamma^B (J) = \sum_{C \in {\cal C}on (J)} <C ,  B> - \sum_{C \in {\cal C}on_{red} (J)} <C ,  B>
\mult_B (C) \in \Z.$$
\begin{theo}
\label{theocon}
The integer $\Gamma^B (J)$ does not depend on the generic choice of $J \in \R {\cal J}_\omega$. Moreover,
it is invariant under isotopy of $B = B_1 \cup \dots \cup B_5$.
\end{theo}
In particular, during such an isotopy, the five curves $L_1 , \dots , L_5$ have to remain transversal 
to each other. Note that there are only finitely many isotopy classes of five real conics in the plane.
How does $\Gamma^B$ depend on the isotopy classes will be studied in \S \ref{subsectdepend}, Proposition 
\ref{propdepend}. The integer $\Gamma^B$ is computed in the following cases.
\begin{prop}
\label{propcomput}
Let $B_1 , \dots , B_5$ be five disjoint disks in $\R P^2$, then $\Gamma^B = 272$. The same holds
when $B_1 , \dots , B_5$ are close to a generic configuration of five real double lines of the plane.
\end{prop}
Here, a disk is said to be close to a double line with equation $y^2 = 0$ in the plane if it has an
equation of the form $\{ y^2 \leq \epsilon^2 x^2 - \delta \}$ for small $\epsilon$ and
$\delta$'s.
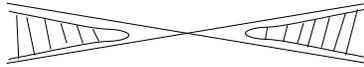
$$\vcenter{\hbox{\begin{picture}(0,0)%
\includegraphics{rel1.pstex}%
\end{picture}%
\setlength{\unitlength}{3315sp}%
\begingroup\makeatletter\ifx\SetFigFont\undefined%
\gdef\SetFigFont#1#2#3#4#5{%
  \reset@font\fontsize{#1}{#2pt}%
  \fontfamily{#3}\fontseries{#4}\fontshape{#5}%
  \selectfont}%
\fi\endgroup%
\begin{picture}(2724,474)(2239,-5023)
\end{picture}
}}$$
\begin{cor}
\label{corcomput}
Let $L_1 , \dots , L_5$ be five real generic plane conics whose isotopy class is given by Proposition
\ref{propcomput}. Then, the number of real conics tangent to $L_1 , \dots , L_5$ is bounded from
below by $32$.
\end{cor}
{\bf Proof:}

The number of lines tangent to two different generic conics is four, they correspond to
the intersection points between
the two dual conics. The number of real reducible conics tangent to four out of the five conics
$L_1 , \dots , L_5$ thus does not exceed $240 = 5*3*4*4$. The result follows now from the definition
of $\Gamma^B$ and Proposition \ref{propcomput}. $\square$\\

Hence, this Corollary \ref{corcomput} provides lower bounds in real enumerative geometry. Note that
this number of real conics does not admit any non trivial upper bound. Indeed,
F. Ronga, A. Tognoli and T. Vust have found a configuration of five real conics close to the double
edges of some pentagon such that all the $3264$ conics tangent to them are real, see \cite{RTV}.

The paper is organized as follows. The first paragraph is devoted to the construction
of the moduli space $\R {\cal M}^d_L$ of real rational pseudo-holomorphic curves which
realize the homology class $d$ and are tangent to $L$. The second paragraph is devoted to the proof
of the results of \S \ref{subsectresults} and the third paragraph to the proof
of the results of \S \ref{subsectcon}.\\

{\bf Acknowledgements:}

This work was initiated during my stay at the Mathematical Sciences Research Institute
in spring $2004$. I would like to acknowledge MSRI for the excellent 
working conditions it provided to me. I also would like to acknowledge Y. Eliashberg, I. Itenberg
and O. Viro for the fruitful discussions we had.

\tableofcontents

\section{Moduli space of real rational pseudo-holomorphic curves tangent to $L$}
\label{sectmoduli}

Let $d \in H_2 (X ; \Z)$ be such that $(c_X)_* d = -d$ and $c_1 (X) d > 1$, $c_1 (X) d \neq 4$. 
Let $\tau$ be an order two permutation of the set $\{ 1 , \dots , c_1 (X) d -2 \}$ having one
fixed point at least. Let
$c_\tau : (x_1 , \dots , x_{c_1 (X) d -2}) \in X^{c_1 (X) d -2} \mapsto
(c_X (x_{\tau (1)}) , \dots , c_X (x_{\tau (c_1 (X) d -2)}) \in X^{c_1 (X) d -2}$
be the associated real structure of $X^{c_1 (X) d -2}$, its fixed point set 
is denoted by $\R_\tau X^{c_1 (X) d -2}$. Let $L \subset \R X$ be a smooth curve such that
$[L] = 0 \in H_1 (\R X ; \Z / 2\Z)$ and $B \subset \R X$ be a surface having $L$ as a
boundary. Finally, let $g$ be a riemannian metric on $X$, invariant under $c_X$ and
for which $L$ is a geodesic. We denote by $\nabla$ the associated Levi-Civita connection on $TX$.

\subsection{Moduli space $\R {\cal P}^*_L$ of real rational pseudo-holomorphic maps tangent to $L$}
\label{subsectmodulitan}

Let $S$ be an oriented sphere of dimension two and $conj$ be a smooth involution conjugated to 
the complex
conjugation of $\C P^1$. Denote by $\R S$ the fixed point set of $conj$ and by
$\R {\cal J}_S$ the space of complex 
structures of class $C^l$ of $S$ which are compatible with its orientation and for which $conj$
is $J$-antiholomorphic. Let $\xi \in \R S$, $\stackrel{\to}{\xi} \in T_\xi \R S \setminus \{ 0 \}$
and $\underline{z} = (z_1 , \dots , z_{c_1 (X) d -2}) \in S^{c_1 (X) d -2}$ be an ordered 
set of $c_1 (X) d -2$ distinct points of $S \setminus \{ \xi \}$. We assume that $\underline{z}$ is
globally invariant under $conj$ and that the permutation of $\{ 1 , \dots , c_1 (X) d -2 \}$
induced by $conj$ is $\tau$. We set
$$\R {\cal P}_L = \{ (u , J_S , J , \underline{x}) \in L^{k,p} (S,X) \times \R {\cal J}_S
\times \R {\cal J}_\omega \times  \R_\tau X^{c_1 (X) d -2} \, | \, u_* [S] = d \, , \,
u(\underline{z}) = \underline{x}$$ $$du + J \circ du \circ J_S = 0 \, , \,
c_X \circ u = u \circ conj \, , \, u(\xi) \in L \text{ and } d_\xi u (\stackrel{\to}{\xi}) \in 
T_{u(\xi)} L \},
$$
where $1 \ll k \ll l$ is large enough and $p > 2$. 

Let $\R {\cal P}^*_L \subset \R {\cal P}_L$ be the
space of {\it non multiple} pseudo-holomorphic maps, that is the space of quadruples
$(u , J_S , J , \underline{x})$ for which $u$ cannot be written $u' \circ \Phi$
where $\Phi : S \to S'$ is a non trivial ramified covering and $u' : S' \to X$
a pseudo-holomorphic map. 

\begin{prop}
\label{propp*}
The space $\R {\cal P}^*_L$ is a separable Banach manifold of class $C^{l-k}$ with
tangent bundle $$T_{(u , J_S , J , \underline{x})} \R {\cal P}^*_L =
\{(v,\dot{J}_S , \dot{J} , \stackrel{.}{\underline{x}}) 
\in  T_{(u , J_S , J , \underline{x})}
(L^{k,p} (S,X) \times \R {\cal J}_S
\times \R {\cal J}_\omega \times  \R_\tau X^{c_1 (X) d -2})
 \, | $$ $$v(\underline{z}) = 
\stackrel{.}{\underline{x}} \, , 
dc_X \circ v = v \circ conj \, , \, v(\xi) \in T_{u(\xi)} L \, , \,
\nabla_{\stackrel{\to}{\xi}} v \in T_{u(\xi)} L \text{ and }
Dv + J \circ du \circ \dot{J}_S +
\dot{J} \circ du \circ J_S = 0 \}.$$
\end{prop}
Here, $T_u L^{k,p} (S,X) = \{ v \in L^{k,p} (S,E_u) \}$ where $E_u = u^* TX$ and
$D : v \in L^{k,p} (S,E_u) \mapsto \nabla v + J \circ \nabla v \circ  J_S +
\nabla_v J \circ du \circ  J_S \in L^{k-1,p} (S, \Lambda^{0,1} S \otimes E_u)$ is
the associated Gromov operator (see \cite{McDSal}, Proposition $3.1.1$).\\

{\bf Proof:}

If we remove the Cauchy-Riemann equation $du + J \circ du \circ J_S = 0$ from the definition
of $\R {\cal P}^*_L$, then the corresponding space $\R {\cal A}^*_L$ is a separable Banach manifold 
of class $C^{l-k}$. After differentiation, the equation $d_\xi u (\stackrel{\to}{\xi}) =
\lambda (u) \zeta (u)$, where $\zeta$ is a unitary vector field tangent to $L$, becomes
$\nabla_{\stackrel{\to}{\xi}} v = d \lambda (u) \zeta (u) + \lambda (u) \nabla_v \zeta (u)$.
Since $L$ is a geodesic for $g$ and $v$ is collinear to $\zeta$, the term $\nabla_v \zeta (u)$
vanishes and $\nabla_{\stackrel{\to}{\xi}} v \in T_{u(\xi)} L$. We have to prove that the space of
non multiple pseudo-holomorphic maps is a Banach submanifold of $\R {\cal A}^*_L$. This follows
from the fact that the section $\sigma_{\overline{\partial}} : (u , J_S , J , \underline{x})
\mapsto du + J \circ du \circ J_S$ of the bundle $L^{k-1 , p} (S , \Lambda^{0,1} S \otimes E_u)$
vanishes transversely, the proof of the latter being the same as the one of Proposition $3.2.1$
of \cite{McDSal}. $\square$

\subsection{Normal sheaf}
\label{subsectionnormal}

Remember that the $\C$-linear part of the Gromov operator $D$ 
is some $\overline{\partial}$-operator denoted by $\overline{\partial}$. The latter induces a
holomorphic structure on the bundle $E_u = u^* TX$ which turns the morphism $du : TS \to
E_u$ into an injective homomorphism of analytic sheaves (see \cite{IShev}, Lemma $1.3.1$).
Likewise, the $\C$-antilinear part of $D$ is some order $0$ operator 
denoted by $R$ and defined by the formula $R_{(u,J_S , J, \underline{x})} (v) = 
N_J (v, du)$ where $N_J$ is the Nijenhuis tensor of $J$. 
Denote by ${\cal O}_S (E_u)$ (resp. ${\cal O}_S (TS)$) the sheaf of analytic $\Z/2\Z$-equivariant
sections of $E_u$ (resp. $TS$). Also, denote by 
${\cal N}_u$ the quotient sheaf ${\cal O}_S (E_u)/du({\cal O}_S (TS))$ so that it fits in
the following exact sequence of analytic sheaves $0 \to {\cal O}_S (TS)\to 
{\cal O}_S (E_u) \to {\cal N}_u \to 0$. Denote by $E_u^L$ the sheaf of $\Z/2\Z$-equivariant
analytic sections of $E_u$ which satisfy $v(\xi) \in T_{u(\xi)} L$ and
$\nabla_{\stackrel{\to}{\xi}} v \in T_{u(\xi)} L$.

\begin{lemma}
\label{lemmadu}
Let $w$ be a real vector field on $S$ which vanishes at $\xi$. Then $du(w) \in E_u^L$.
\end{lemma}
{\bf Proof:}

Denote by $v = du(w)$, then $v \in E_u$ and $v (\xi) = du (w (\xi)) = 0$. Moreover,
$\nabla_{\stackrel{\to}{\xi}} v = (\nabla_{\stackrel{\to}{\xi}} du) (w) + 
du (\nabla_{\stackrel{\to}{\xi}} w)$. The first term vanishes since $w (\xi) = 0$ and the second
belongs to $T_{u(\xi)} L$. $\square$\\

We denote by ${\cal N}_u^L$ the quotient sheaf ${\cal O}_S (E_u^L) /
du ({\cal O}_S (TS_{-\xi}))$, so that we have the exact sequence
$0 \to {\cal O}_S (TS_{-\xi}) \to {\cal O}_S (E_u^L) \to {\cal N}_u^L \to 0$.
Denote by $\widetilde{E}_u^L$ the sheaf of $\Z/2\Z$-equivariant
sections of $E_u$ for which $v (\xi) \in T_{u(\xi)} L$, so that 
$E_u^L \subset \widetilde{E}_u^L$.

\begin{lemma}
\label{lemmaetilde}
If $d_\xi u \neq 0$, then the quotient of $\widetilde{E}_u^L$ by $du ({\cal O}_S (TS))$ is
the sheaf ${\cal N}_{u, -\xi}$. If $d_\xi u = 0$, but $\nabla_{\stackrel{\to}{\xi}} du
(\stackrel{\to}{\xi}) \notin T_{u(\xi)} L$, then this quotient is the sheaf
$N_u = {\cal O}_S (E_u) / du ({\cal O}_S (TS) \otimes {\cal O}_S (\xi))$.
\end{lemma}
{\bf Proof:}

The first part follows from the fact that the condition $v (\xi) \in T_{u(\xi)} L$ for a section
$v$ of $E_u$ reads in the quotient as a section of ${\cal N}_u$ which vanishes at $\xi$,
since $T_{u(\xi)} L \subset Im (d_\xi u)$. In the second case, the normal sheaf ${\cal N}_u$
splits as $N_u \oplus {\cal N}_u^{sing}$, where ${\cal N}_u^{sing}$ is the skyscrapper
part $du({\cal O}_S (TS) \otimes {\cal O}_S (\xi)) / du({\cal O}_S (TS))$. From the hypothesis,
the cuspidal point at $\xi$ is non-degenerated and has a tangent line distinct from $T_{u(\xi)} L$.
Thus, $du(TS \otimes {\cal O}_S (\xi)) \not\subset \widetilde{E}_u^L$ and the
skyscrapper part ${\cal N}_u^{sing}$ does not belong to the quotient $\widetilde{E}_u^L /
du ({\cal O}_S (TS))$. The projection $\widetilde{E}_u^L /
du ({\cal O}_S (TS)) \subset {\cal N}_u$ onto $N_u$ induced by $N_u \oplus {\cal N}_u^{sing} \to N_u$
provides the required isomorphism. $\square$\\

As soon as $d_\xi u \neq 0$, we deduce the exact sequence $0 \to {\cal O}_S (TS_{-\xi}) \to 
\widetilde{E}_u^L \to {\cal N}_{u,-\xi} \oplus T_\xi L \to 0$, where $T_\xi L$ is the skyscrapper
sheaf $du ({\cal O}_S (TS)) / du ({\cal O}_S (TS_{-\xi}))$. We deduce the inclusion
${\cal N}_u^L \subset {\cal N}_{u,-\xi} \oplus T_\xi L$.

\begin{prop}
\label{propprepcrit}
1) Assume that $d_\xi u \neq 0$. Then, the skyscrapper
part $T_\xi L$ is included in ${\cal N}_u^L$ if and only if $\nabla_{\stackrel{\to}{\xi}} du 
(\stackrel{\to}{\xi}) \in T_{u(\xi)} L$, that is if $u(\xi)$ is a degenerated point of contact
between $u(S)$ and $L$. In this case, the projection ${\cal N}_{u,-\xi} \oplus T_\xi L \to
{\cal N}_{u,-\xi}$ restricted to ${\cal N}_u^L$ has image ${\cal N}_{u,-2 \xi}$. Otherwise,
this projection establishes an isomorphism between ${\cal N}_u^L$ and ${\cal N}_{u,-\xi}$.

2) Assume that $d_\xi u = 0$ but $\nabla_{\stackrel{\to}{\xi}} du (\stackrel{\to}{\xi})
\notin T_{u(\xi)} L$. Then, the sheaf ${\cal N}_u^L$ is isomorphic to ${\cal N}_{u}$.
\end{prop}
{\bf Proof:}

If $d_\xi u \neq 0$, the skyscrapper
part $T_\xi L$ is generated by $du ({\cal O}_S (TS))$. Let $w$ be a real vector field on $S$, we
have to see under which condition $du (w) \in E_u^L$. From the relation
\begin{eqnarray}
\label{equation1}
\nabla_{\stackrel{\to}{\xi}} (du (w)) & = & (\nabla_{\stackrel{\to}{\xi}} du) (w) + 
du (\nabla_{\stackrel{\to}{\xi}} w),
\end{eqnarray}
it is necessary and sufficient that $\nabla_{\stackrel{\to}{\xi}} du (\stackrel{\to}{\xi})
\in T_{u(\xi)} L$. In this case, the connection $\nabla$ induces at $\xi$ a derivation
$\nabla^{\xi}$ of sections of the sheaf ${\cal N}_{u}$ such that the
relations $v(\xi) \in T_{u(\xi)} L$ and
$\nabla_{\stackrel{\to}{\xi}} v \in T_{u(\xi)} L$ reads in the quotient
$v(\xi) = 0$ and $\nabla^{\xi} v = 0$. Thus, the projection ${\cal N}_{u,-\xi} \oplus T_\xi L \to
{\cal N}_{u,-\xi}$ restricted to ${\cal N}_u^L$ has image ${\cal N}_{u,-2 \xi}$. Otherwise,
it induces an isomorphism.

Assume now that $d_\xi u = 0$ but $\nabla_{\stackrel{\to}{\xi}} du (\stackrel{\to}{\xi})
\notin T_{u(\xi)} L$. Then, if $w$ is a real vector field on $S$ such that $w (\xi) \neq 0$,
$du (w) \notin E_u^L$ from (\ref{equation1}). The exact sequence $0 \to {\cal O}_S (TS) \to 
\widetilde{E}_u^L \to {\cal N}_{u} \to 0$ restricts thus as
$0 \to {\cal O}_S (TS_{-\xi}) \to 
E_u^L \to {\cal N}_{u} \to 0$, hence the result. $\square$\\

Denote by ${\cal O}_S (TS_{- \underline{z}})$ (resp. ${\cal O}_S (E_{u , 
- \underline{z}})$, ${\cal O}_S (E^L_{u , 
- \underline{z}})$, ${\cal N}_{u , 
- \underline{z}}$, ${\cal N}^L_{u , 
- \underline{z}}$,  $T_{u(\xi) , - \underline{z}} L$) the subsheaf of sections of 
${\cal O}_S (TS)$ 
(resp. ${\cal O}_S (E_{u})$, ${\cal O}_S (E^L_{u})$, ${\cal N}_{u}$, 
${\cal N}^L_{u}$,  $T_{u(\xi)} L$) which vanish at 
$\underline{z}$.
Remember that the operator $D : L^{k,p} (S , E^L_{u , -\underline{z}}) \to
L^{k-1,p} (S , \Lambda^{0,1} S \otimes E_{u})$ induces a quotient
operator $\overline{D} : L^{k,p} (S , {\cal N}^L_{u , -\underline{z}}) :=
L^{k,p} (S , E^L_{u , -\underline{z}}) / du(L^{k,p} (S , TS_{-\underline{z}})
\to$ $L^{k-1,p} (S , \Lambda^{0,1} S \otimes N^L_{u})$.
From the short exact sequence of complexes
$$\begin{array}{ccccccccc}
0 & \to & L^{k,p} (S , TS_{-\xi -\underline{z}}) & \stackrel{du}{\to} & 
L^{k,p} (S , E^L_{u , -\underline{z}}) & \to & 
L^{k,p} (S , {\cal N}^L_{u , -\underline{z}})
& \to & 0 \\
&&\downarrow \overline{\partial}_S && \downarrow D && \downarrow \overline{D}&&\\
0 & \to & L^{k-1,p} (S , \Lambda^{0,1} S \otimes TS) & \stackrel{du}{\to} & 
L^{k-1,p} (S , \Lambda^{0,1} S \otimes E^L_{u}) & \to & 
L^{k-1,p} (S , \Lambda^{0,1} S \otimes N^L_{u})
& \to & 0, \\
\end{array}$$
we deduce the long exact sequence
$0 \to H^0 (S , TS_{-\xi -\underline{z}}) \to H^0_D (S , E^L_{u , -\underline{z}}) \to 
H^0_{\overline{D}} (S , {\cal N}^L_{u , -\underline{z}}) \to
H^1 (S , TS_{-\xi -\underline{z}}) \to H^1_D (S , E^L_{u , -\underline{z}}) \to 
H^1_{\overline{D}} (S , {\cal N}^L_{u , -\underline{z}}) \to 0$,
where $H^0_D$, $H^0_{\overline{D}}$ (resp. $H^1_D$, $H^1_{\overline{D}}$) denote
the kernels (resp. cokernels) of the operators $D$, $\overline{D}$ on the
associated sheaves. In particular,
\begin{eqnarray*}
\ind_\R (\overline{D}) & = & \ind_\R (D) - \ind_\R (\overline{\partial}_S) \\
 & = & ( c_1 (X)d + 2 - 2 - 2 \# \underline{z}) - (3 - 1 - \# \underline{z}) \\
& = & 0.
\end{eqnarray*}

\subsection{Moduli space of real rational pseudo-holomorphic curves tangent to $L$}

Denote by ${\cal D}iff_\R^+ (S,z, \xi)$ the group of diffeomorphisms of class $C^{l+1}$ of $S$,
which preserve the orientation, fix $\underline{z} \cup \{ \xi \}$ and commute with $conj$.
This group acts on $\R {\cal P}^*_L $ by
$$\phi . (u,J_S , J , \underline{x})  = 
(u \circ \phi^{-1} , (\phi^{-1})^*J_S , J , \underline{x}),$$
where $(\phi^{-1})^* J_S = d\phi \circ J_S \circ d\phi^{-1}$.
Denote by $\R {\cal M}^d_L$ the quotient of 
$\R {\cal P}^*_L$  by this action.
The projection $\pi : (u,J_S , J , \underline{x}) \in \R {\cal P}^*_L  \mapsto (J,\underline{x}) 
\in {\cal J}_\omega \times X^{c_1 (X)d - 2}$
induces on the quotient a projection $\R {\cal M}^d_L \to \R {\cal J}_\omega \times \R_\tau
X^{c_1 (X)d - 2}$ still denoted by $\pi$. 
\begin{prop}
\label{proppi}
The space  $\R {\cal M}^d_L$ is a separable Banach manifold
of class $C^{l-k}$, and $\pi$ is Fredholm of vanishing index. Moreover, if
$[u,J_S , J , \underline{x}] \in \R {\cal M}^d_L $, then we have the isomorphisms
$\ker d\pi|_{(u,J_S , J , \underline{x})} \cong H^0_D (S , {\cal N}^L_{u , -\underline{z}})$
and $\coker d\pi|_{(u,J_S , J , \underline{x})} \cong 
H^1_D (S , {\cal N}^L_{u , -\underline{z}})$.
\end{prop}

{\bf Proof:}

The proof is analogous to the one of Corollary $2.2.3$ of \cite{Shev} and Proposition $3.2.1$
of \cite{McDSal}.
The action of ${\cal D}iff_\R^+ (S,z)$ on $\R {\cal P}^*_L$ 
is smooth, fixed point free and admits a closed supplement. From Proposition \ref{propp*}
thus follows that $\R {\cal M}^d_L$ is a separable Banach manifold of class $C^{l-k}$. Moreover, 
\begin{eqnarray*}
\ker d\pi|_{[u,J_S , J , \underline{x}]} & = &
\{ (v, \dot{J}_S , 0 , 0 ) \in T_{(u,J_S , J , \underline{x})} \R {\cal P}^*_L
\, | \, v(\underline{z}) = 0 \}/T_{Id} {\cal D}iff_\R^+ (S,z, \xi) \\
 & = & \{ v \in L^{k , p} (S , E^L_{u , -\underline{z}}) \, | \, \exists \phi \in 
L^{k-1,p} (S , \Lambda^{0,1} S \otimes TS) \, , \, D v = du (\phi) \}/ du (L^{k,p} (S , 
TS_{-\underline{z}})) \\
 & = & H^0_D (S , {\cal N}^L_{u , -\underline{z}}).
\end{eqnarray*}
Likewise, 
\begin{eqnarray*}
Im d\pi|_{[u,J_S , J , \underline{x}]} & = & \{ (\dot{J} , \stackrel{.}{\underline{x}})
\in T_J \R {\cal J}_\omega \times T_{\underline{x}} \R_\tau X^{c_1 (X)d - 2} \, | \, \exists 
(v , \dot{J}_S) \in  L^{k , p} (S , E^L_{u , -\underline{z}}) \times T_{J_S} \R {\cal J}_S,\\
&&  Dv + J \circ du  \circ \dot{J}_S 
= - \dot{J} \circ du \circ J_S ,
\,  v(\underline{z}) = \stackrel{.}{\underline{x}} \}, \text{ so that } \\
\coker d\pi|_{[u,J_S , J , \underline{x}]} &\cong&  L^{k -1, p} (S , \Lambda^{0,1} S \otimes E^L_u)  
\times T_{\underline{x}} 
\R_\tau X^{c_1 (X)d -2} /
Im(\widehat{D} \times ev),
\end{eqnarray*}
where $\widehat{D} : (v,\dot{J}_S) \in L^{k , p} (S , E^L_u)
\times T_{J_S} \R {\cal J}_S \mapsto Dv + J \circ du \circ \dot{J}_S \in 
L^{k -1, p} (S , \Lambda^{0,1} S \otimes E^L_u)$ and $ev : v \in L^{k , p} (S , E^L_u)
\mapsto v(\underline{z}) \in T_{\underline{x}} \R_\tau X^{c_1 (X)d -2}$.
In particular, $Im d\pi|_{[u,J_S , J , \underline{x}]}$ is closed and $\pi$ is Fredholm.
By definition, $\coker D = H^1_D (S , E^L_u)$. From the short exact sequence
$0 \to E^L_{u , - \underline{z}} \to E^L_u \stackrel{ev}{\to} T_{\underline{x}} 
\R_\tau X^{c_1 (X)d -2} \to 0$, we deduce the long exact sequence $\to H^0_D (S , E^L_u) \to
H^0 (S ,  T_{\underline{x}} 
\R_\tau X^{c_1 (X)d -2}) \to H^1_D (S , E^L_{u , - \underline{z}}) \to  
H^1_D (S , E^L_u) \to 0$.
Hence, the cokernel of $D \times ev$ in $T_J \R {\cal J}_\omega \times
T_{\underline{x}} \R_\tau X^{c_1 (X)d -2}$ is isomorphic to $H^1_D (S , E^L_{u , - \underline{z}})$.
From the long exact sequence given at the end of \S \ref{subsectionnormal}, we deduce that
the cokernel of $\widehat{D} \times ev$ and hence the one of $d\pi|_{[u,J_S , J , \underline{x}]}$
is isomorphic to $H^1_D (S , {\cal N}^L_{u , -\underline{z}})$. $\square$
\begin{cor}
\label{corcrit}
The critical points $[u,J_S , J , \underline{x}]$ of $\pi$ are those for which $u(S)$ has
a point of contact of order greater than two with $L$ at $u(\xi)$ or $u$ has a cuspidal point 
outside $\xi$. $\square$
\end{cor}

\subsection{Generic critical points of $\pi$ are non degenerated}
\label{subsectcritical}

\begin{theo}
\label{theocrit}
Let $[u,J_S , J , \underline{x}] \in \R {\cal M}^d_L$ be such that $u(S)$ has
a point of contact of order two with $L$ at $u(\xi)$  and a unique real ordinary cuspidal point 
outside $\xi$. Then, $[u,J_S , J , \underline{x}]$ is a non degenerated
critical points of $\pi$. The same holds if $u(S)$ is immersed but has a point of contact of order three
with $L$ at $u(\xi)$.
\end{theo}
The critical points of $\pi$ which appear in this Theorem \ref{theocrit} are
said to be {\it generic}.\\

{\bf Proof:}

The proof of the first part of this theorem is the same as the one of Lemma $2.13$ of
\cite{Wels1}, it is not reproduced here. Let $[u,J_S , J , \underline{x}] 
\in \R {\cal M}^d_L$ be such that $u(S)$ is immersed but has a point of contact of order three
with $L$ at $u(\xi)$. We have to prove that the quadratic form $\nabla 
d\pi |_{[u,J_S , J , \underline{x}]} : \ker d\pi |_{[u,J_S , J , \underline{x}]}
\times \ker d\pi |_{[u,J_S , J , \underline{x}]} \to \coker 
\ker d\pi |_{[u,J_S , J , \underline{x}]}$ is non degenerated.
We saw in the proof of Proposition \ref{proppi} that the kernel and cokernel
of the map $d\pi$ are the same as the ones of the morphism $-\widehat{D}_\R :
(v , \dot{J}_S , \dot{J} , \stackrel{.}{x}) \in T_{[u,J_S , J, \underline{x}]}
\R {\cal M}^d_L   \mapsto \dot{J} \circ du \circ J_S \in 
L^{k-1 , p} (S , \Lambda^{0,1} S \otimes N_{u}^L)$. From the relation 
$Dv + J \circ du \circ 
\dot{J}_S + \dot{J} \circ du \circ J_S = 0$, we deduce that
$\widehat{D}_\R (v , \dot{J}_S , \dot{J} , \stackrel{.}{x}) = Dv +
J \circ du \circ \dot{J}_S$. We then have to prove that
$\nabla \widehat{D}_\R|_{[u,J_S , J , \underline{x}]} : 
H^0_D (S , {\cal N}_{u, -\underline{z}}^L)^2 \to
H^1_D (S , {\cal N}_{u, -\underline{z}}^L)$ is non degenerated. Let 
$(v, \dot{J}_S , 0 , 0)$ be a generator of $H^0_D (S , {\cal N}_{u, -\underline{z}}^L)$. From
Propositions
\ref{propprepcrit} and \ref{proppi}, $v = du(w)$ for some real vector field $w$  
on $S$ which does not vanish at $\xi$. We can assume that $\dot{J}_S$ vanishes in a neighbourhood of
$z \cup \xi$. After differentiation of the relation 
$D \circ du = du \circ \overline{\partial}_S$,
we deduce
$$\nabla_v D \circ du + D \circ
(\nabla_v du) +  \nabla_{\dot{J}_S} D \circ du = (\nabla_v du) \circ \overline{\partial}_S 
 \mod(Im(du)).$$
Moreover, $\nabla_{(v, \dot{J}_S , 0 , 0)} \widehat{D} = 
\nabla_v D + (\nabla_v du) \circ J_S \circ \dot{J}_S + \nabla_{\dot{J}_S} D
\mod(Im(du))$.
Since the relation $Dv + J \circ du \circ \dot{J}_S =0$ forces
$\overline{\partial}_S (w) + J_S \dot{J}_S =0$, we get (compare Lemma $2.13$ of \cite{Wels1}
and Theorem $1.8$ of \cite{Wels2})
$$(\nabla_{(v , 
\dot{J}_S , 0 , 0)} \widehat{D} )(v) + D (\nabla_{v} du) (w) = 0 
\mod(Im(du)).$$
From Proposition \ref{propprepcrit}, ${\cal N}_{u , -\underline{z}}^L \cong {\cal N}_{u,- \underline{z}
-2 \xi}$. From Riemann-Roch duality, $H^1_D (S , {\cal N}_{u,- \underline{z}
-2 \xi})_{+1}^* \cong H^0_D (S , K_S \otimes {\cal N}_{u,- \underline{z}
-2 \xi})_{-1}$, see \cite{Wels1}, Lemma $1.7$. Let $\psi$ be a generator of
$H^0_D (S , {\cal N}_{u,- \underline{z} -2 \xi})_{-1}$ so that 
$D^* \psi$ is a linear combination of Dirac sections of $N_u^*$ at $z \cup \xi$ as well as of the
derivative $\delta'_{\xi}$ of the Dirac section at $\xi$. Note that since $H^0_D (S , 
K_S \otimes {\cal N}_{u,- \underline{z} - \xi}) = 0$, the coefficient $a_{\xi}$ of $\delta'_{\xi}$
in $D^* \psi$ does not vanish.
We have:
\begin{eqnarray*}
<\psi , \nabla d\pi ((v, \dot{J}_S) , (v, \dot{J}_S))> &=&
-<\psi , \nabla \widehat{D} ((v, \dot{J}_S) , (v, \dot{J}_S))>\\
&=& <\psi , D (\nabla_{v} du) (w) >\\
&=& <D^* \psi , (\nabla_{v} du) (w) >
\end{eqnarray*}
Choose a local chart at $u(\xi)$ such that $L$ is conjugated to the first coordinate axis of $\R^2
\subset \C^2$. Without loss of generality, we can assume that the first coordinate axis is
$J$-holomorphic and that the metric $g$ is constant in this chart, so that $\nabla = d$.
The map $u$ writes then $z \mapsto ((z-\xi) +  o(|z-\xi|) , (z-\xi)^3 +  o(|z-\xi|^3))$
in a neighbourhood of $\xi$. Thus, $\nabla_v du (w) = \nabla_w v$, considered as
a section of the normal bundle of $u$, has a simple zero at $\xi$. Since $w$ vanishes at $z$,
we deduce that $<D^* \psi , (\nabla_{v} du) (w) > = a_{\xi} < \delta'_{\xi} , (\nabla_{v} du) (w) >$.
Now since the vanishing order of $\nabla_v du (w)$ at $\xi$ is one, 
$< \delta'_{\xi} , (\nabla_{v} du) (w) > \neq 0 $, hence the result. $\square$

\subsection{Gromov compactification $\R \overline{\cal M}^d_L$ of $\R {\cal M}^d_L$}
\label{subsectcompact}

The projection $\pi : \R {\cal M}^d_L \to \R {\cal J}_\omega \times 
\R_\tau X^{c_1 (X)d - 2}$ is not proper in general. Its lack of properness is described by
the following lemma which follows from Gromov's compactness Theorem (see \cite{McDSal}, Theorem
5.5.5).
\begin{lemma}
\label{lemmacompact}
Let $[u^n,J_S^n , J^n , \underline{x}^n]$ be a sequence of elements of $\R {\cal M}^d_L$ such
that $(J^n , \underline{x}^n)$ converges to $(J^\infty , \underline{x}^\infty)$. Then, after
possibly extracting a subsequence, we have one of the following:

1) This sequence $[u^n,J_S^n , J^n , \underline{x}^n]$ converges in $\R {\cal M}^d_L$.

2) The sequence $u^n (S)$ converges to some irreducible curve, tangent to $L$, but
the point of contact belongs to $\underline{x}^\infty$.

3) The sequence $u^n (S)$ converges to some reducible curve. Moreover, in this case, the reducible curve
is either tangent to $L$, or has two of its irreducible components which intersect on $L$. $\square$
\end{lemma}

\section{Proofs of Theorems \ref{maintheo} and \ref{theorelation}}
\label{sectionproofs}

Let $(J^0 , \underline{x}^0)$ and $(J^1 , \underline{x}^1)$ be two generic elements of
$\R {\cal J}_\omega \times \R_\tau X^{c_1 (X)d - 2}$ so that the integers $\Gamma_r^{d,B}
(J^0 , \underline{x}^0)$ and $\Gamma_r^{d,B} (J^1 , \underline{x}^1)$ are well defined.
We have to prove that they coincide.

\subsection{Choice of a path $\gamma$}
\label{subsectchoice}

Remember that by definition, a {\it stratum of codimension $k \geq 0$} of a separable Banach manifold
$M$ is the image of a separable Banach manifold $L$ under a Fredholm map $\Phi$ of index $-k$
such that the limits of sequences $\Phi (x_n)$ where $(x_n)_{n \in \N}$ diverges in $N$ belong
to a countable union of strata of higher codimensions. In particular, $\Phi$ is not assumed to
be proper.

\begin{prop}
\label{propgeneric}
The subset of elements $[u, J_S , J , \underline{x}]$ of $\R {\cal M}^d_L$ for which $u(S)$
has only transversal double points as singularities, outside $\underline{x} \cup L$, and a unique
point of contact of order two with $L$, is a dense open subset of $\R {\cal M}^d_L$.
The four followings are substrata of codimension one of $\R {\cal M}^d_L$.

1) Curves having  only transversal double points as singularities, outside $\underline{x} \cup L$, 
and a unique point of contact with $L$ which is of order three.

2) Curves having a unique real ordinary cusp and transversal double points as singularities, 
outside $\underline{x} \cup L$, and a unique
point of contact with $L$ which is of order two.  

3) Curves having a unique real ordinary cusp on $L$ and transversal double points outside 
$\underline{x} \cup L$ as singularities. These curves are not tangent to $L$ and the tangent line of
the curve at the cusp is distinct from the one of $L$.

4) Curves having a real ordinary triple point or real ordinary tacnode or a transversal double point
on $\underline{x} \cup L$ or two points of contact with $L$.

The set of curves not listed above belongs to a countable union of strata of codimension greater than
one of $\R {\cal M}^d_L$.
\end{prop}
{\bf Proof:}

The proof is the same as the one of Proposition $2.7$ of \cite{Wels1}. It is left to the reader. 
$\square$\\

Let $\gamma : t \in [0,1] \mapsto (J^t , \underline{x}^t) \in \R {\cal J}_\omega \times 
\R_\tau X^{c_1 (X)d - 2}$ be a generic path transversal to  $\pi_\R$. Denote by $\R {\cal M}_\gamma
= \R {\cal M}^d_L \times_\gamma [0,1]$, $\R \overline{\cal M}_\gamma$ its Gromov compactification
and $\pi_\gamma : \R \overline{\cal M}_\gamma \to [0,1]$ the associated projection.

\begin{prop}
\label{propcompact}
As soon as $\gamma$ is generic enough, the elements of $\R \overline{\cal M}_\gamma \setminus
\R {\cal M}_\gamma$ are either irreducible curves $[u^t, J_S^t , J^t , \underline{x}^t]$ such
that $\underline{x}^t \cap L$ is non empty, or reducible curves $C^t$ having two irreducible components
$C_1^t$, $C_2^t$,
both real, and only transversal double points as singularities, outside $\underline{x}$. Moreover,
we have the following alternative:

1) Either $C^t$ has a unique
point of contact with $L$ which is of order two and outside its singular points.

2) Or $C^t$ has a unique double point on $L$ which is an intersection point of
$\R C_1^t$ and $\R C_2^t$. In this case, it is not tangent to $L$.

Finally, if we denote by $m_i = \# (\underline{x}^t \cap C_i^t)$ and $d_i = [C_i^t] \in H_2 (X ; \Z)$,
$i \in \{ 1,2 \}$, so that $m_1 + m_2 = c_1 (X) d - 2$, then either $m_1 = c_1 (X) d_1 - 1$ or
$m_1 = c_1 (X) d_1 - 2$.
\end{prop} 
{\bf Proof:} 

The proof is the same as the ones of Proposition $2.9$, Corollary $2.10$ and
Proposition $2.11$ of \cite{Wels1}, as well as Corollary $1.12$ of \cite{Wels2}. It is not
reproduced here. $\square$

\begin{rem}
\label{rem0}
Remember that to cover the case $r = (0, \dots , 0)$, one should take into account real reducible
curves made of two complex conjugated components, see Remark $1.9$ of \cite{Wels2}.
It would then be possible to extend Theorem \ref{maintheo} to this case provided an analog
of Theorem $3.2$ of \cite{Wels2} is proved, see Remark $3.5$ of \cite{Wels2}.
\end{rem}

{\bf From now on}, we fix a choice of $\gamma$ generic enough so that $\R \overline{\cal M}_\gamma$ 
consists of curves listed in Propositions \ref{propgeneric} and \ref{propcompact}.

\subsection{Neighbourhood of curves having an order three point of contact with $L$}
\label{subsectdegtan}

\begin{prop}
\label{propinflectif}
Let $C=[u, J_S , J , \underline{x}] \in \R {\cal M}_\gamma$ be a curve having an order three point 
of contact with $L$ and $t_0 = \pi_\gamma (C)$. Then, there exist $\eta > 0$ and a neighbourhood 
$W$ of $C$ in $\R {\cal M}_\gamma$ such that for every $t \in ]t_0 - \eta , t_0[$,
$\pi_\gamma^{-1} (t) \cap W$ consists of two curves $C_t^+$, $C_t^-$ having same mass and for
which $<C_t^+ , B> = - <C_t^- , B>$ and for every $t \in ]t_0 , t_0 + \eta [$,
$\pi_\gamma^{-1} (t) \cap W = \emptyset$, or vive versa.
\end{prop} 
$$\vcenter{\hbox{\input{rel2.pstex_t}}}$$
{\bf Proof:}

From Theorem \ref{theocrit}, $C$ is a non degenerated critical point of $\pi_\gamma$. Since
$\R {\cal M}_\gamma$ is of dimension one, this implies that there exist $\eta > 0$ and a neighbourhood 
$W$ of $C$ in $\R {\cal M}_\gamma$ such that for every $t \in ]t_0 - \eta , t_0[$,
$\pi_\gamma^{-1} (t) \cap W$ consists of two curves and for every $t \in ]t_0 , t_0 + \eta [$,
$\pi_\gamma^{-1} (t) \cap W = \emptyset$, or vive versa. The only thing to prove is that  in the 
first case, the two curves $C_t^+$, $C_t^-$ have the same mass and satisfy $<C_t^+ , B> = - <C_t^- , B>$.
The former is obvious. Choose a parameterization $\lambda \in ] -\sqrt{\eta} , \sqrt{\eta} [
\mapsto C_\lambda = [u^\lambda , J_S^\lambda , J^\lambda , \underline{x}^\lambda] 
\in \R {\cal M}_\gamma$ such that $\pi_\gamma (C_\lambda) = t_0 - \lambda^2$. Fix a local chart
$0 \in ]-1 , 1[$ of $\xi \in \R S$ and $0 \in \R^2$ of $u^0 (\xi) \in \R X$. We can assume that
in this second chart, $L$ is identified with the first coordinate axis and $B$ with the upper half
plane of $\R^2$. The one parameter family $(u^\lambda)_{\lambda \in ] -\sqrt{\eta} , \sqrt{\eta} [}$
reads as a map $f : (\lambda , z) \in ] -\sqrt{\eta} , \sqrt{\eta} [ \times
]-1 , 1[ \mapsto f(\lambda , z) \in \R^2$. Denote by $f_1 (\lambda , z)$ and  $f_2 (\lambda , z)$
the two coordinates of  $f (\lambda , z)$. These maps of class $C^{l-k}$, satisfy
$f_1 (0, z) = z + o(|z|)$, $f_2 (0, z) = z^3 + o(|z|^3)$, $f_2 (\lambda , 0) = 0$ and
$\frac{\partial}{\partial z} f_2 (\lambda , z)|_{z=0} = 0$. Moreover, $\frac{\partial}{\partial
\lambda} C_\lambda|_{\lambda =0}$ generates the kernel of $d \pi_\gamma|_{C_\lambda}$. It thus follows
from Proposition \ref{proppi} that $\frac{\partial}{\partial \lambda} f(\lambda , z)|_{\lambda =0}
= \frac{\partial}{\partial z} f(\lambda , z)|_{\lambda =0} = (1 + o(1) , 3z^2 + o(|z|^2))$. We
deduce that the order three jet  of $f_2$ writes $f_2 (\lambda , z) = z^2 (z + a\lambda ) +
o(||(\lambda , z)||^3)$, for some $a \in \R^*$. Hence, when $\lambda > 0$ (resp. $\lambda < 0$),
the sign of $f_2 (\lambda , z)$ in a neighbourhood of $z =0$ is the one of $a$ (resp. its opposite).
In particular, as soon as $\lambda \neq 0$, $<C_\lambda , B> = - <C_{-\lambda} , B>$. $\square$

\subsection{Neighbourhood of curves having a cuspidal point}
\label{subsectcusp}

\begin{prop}
\label{propcusp1}
Let $C=[u, J_S , J , \underline{x}] \in \R {\cal M}_\gamma$ be a curve having a real ordinary cusp
outside $L$ and $t_1 = \pi_\gamma (C)$. Then, there exist $\eta > 0$ and a neighbourhood 
$W$ of $C$ in $\R {\cal M}_\gamma$ such that for every $t \in ]t_1 - \eta , t_1[$,
$\pi_\gamma^{-1} (t) \cap W$ consists of two curves $C_t^+$, $C_t^-$ such that $m(C_t^+) = m(C_t^-)
+ 1$ and $<C_t^+ , B> = <C_t^- , B>$ and for every $t \in ]t_1 , t_1 + \eta [$,
$\pi_\gamma^{-1} (t) \cap W = \emptyset$, or vive versa.
\end{prop} 
{\bf Proof:}

From Theorem \ref{theocrit}, $C$ is a non degenerated critical point of $\pi_\gamma$. Since
$\R {\cal M}_\gamma$ is of dimension one, this implies that there exist $\eta > 0$ and a neighbourhood 
$W$ of $C$ in $\R {\cal M}_\gamma$ such that for every $t \in ]t_1 - \eta , t_1[$,
$\pi_\gamma^{-1} (t) \cap W$ consists of two curves and for every $t \in ]t_1 , t_1 + \eta [$,
$\pi_\gamma^{-1} (t) \cap W = \emptyset$, or vive versa. The only thing to prove is that
$m(C_t^+) = m(C_t^-) + 1$. The proof of this is readily the same as the one of Proposition $2.16$
of \cite{Wels1}. It is not reproduced here. $\square$

\begin{prop}
\label{propcusp2}
Let $C=[u, J_S , J , \underline{x}] \in \R {\cal M}_\gamma$ be a curve having a real ordinary cusp
on $L$ and $t_2 = \pi_\gamma (C)$. Then, there exist $\eta > 0$ and a neighbourhood 
$W$ of $C$ in $\R {\cal M}_\gamma$ such that for every $t \in ]t_2 - \eta , t_2 + \eta[ \setminus
\{ t_2 \}$, $\pi_\gamma^{-1} (t) \cap W$ is reduced to one element $\{ C_t \}$. Moreover, 
$<C_t , B>$ does not depend on $t \in ]t_2 - \eta , t_2 + \eta[ \setminus \{ t_2 \}$. Likewise,
$C$ extends to a one parameter family $C_t^{cusp}$ of cuspidal real rational $J^t$-holomorphic curves
which pass through $\underline{x}^t$ and realize $d$. Assume that for $t \in ]t_2 - \eta , t_2 [$
(resp. $t \in ]t_2 , t_2 + \eta [$), $\R C_t^{cusp}$ does not intersect locally $L$ (resp. intersects $L$
locally in two points) near the cusp of $C$. Then for $t \in ]t_2 - \eta , t_2 [$, $m(C_t) = m(C)$ and
for $t \in ]t_2 , t_2 + \eta [$, $m(C_t) = m(C) + 1$.
\end{prop} 
Note that after changing the parameterization $t \mapsto 2t_2 -t$ if necessary, we can always
assume that for $t \in ]t_2 - \eta , t_2 [$
(resp. $t \in ]t_2 , t_2 + \eta [$), $\R C_t^{cusp}$ does not intersect locally $L$ (resp. intersects 
locally $L$ in two points) near the cusp of $C$.
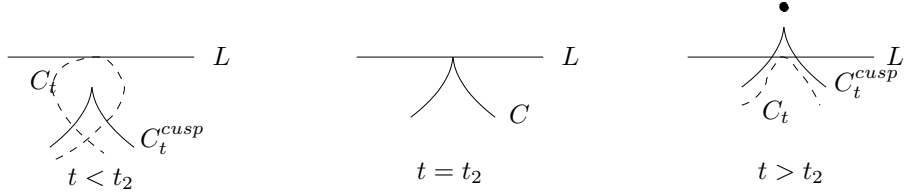
$$\vcenter{\hbox{\input{rel3.pstex_t}}}$$
{\bf Proof:}

Remember that the choice of $\gamma$ implies that the tangent line of $C$ at the cusp is distinct from
the one of $L$.
Without loss of generality, we can assume that $J$, $\underline{x}$ are constant and that $L$
(and the metric $g$) moves along a one parameter family $L_t$ which crosses the cuspidal point of $C$.
This indeed can be realized equivalently by fixing $L$ and having $J$, $\underline{x}$ moving
 along one parameter families $\phi_t^* J$, $\phi_t (\underline{x})$ where $\phi_t$ is some
$\Z/2\Z$-equivariant hamiltonian flow of $X$. The family of curves $C_t^{cusp}$ is then nothing but
the constant family $C$. Moreover, from Proposition $2.16$ of \cite{Wels1}, the curve $C$ extends to
a one parameter family $C^\lambda$, $\lambda \in ]-\epsilon , \epsilon [$, of real rational 
$J$-holomorphic curves which pass through $\underline{x}$ and realize $d$. These curves $C^\lambda$
have an isolated real double point near the cusp of $C$ when $\lambda < 0$ and a non isolated one
when $\lambda > 0$. Moreover, the latter form a one parameter family of loops which fill some
disk of $\R X$ centered at the cusp of $C$ (compare \cite{Wels2}, Lemma $3.3$). This follows from the
fact that the intersection points between two curves in this family are located near their double 
points and at $\underline{x}$. Since for $t \in ]t_2 - \eta , t_2 [$, $L_t$ is locally disjoint from
$C$, there does exist some curve $C^\lambda$, $\lambda > 0$, in this family which is tangent to $L_t$,
as soon as $\eta$ is small enough. It has the same mass as $C$. From Corollary \ref{corcrit},
$C$ is a regular point of $\pi_\gamma$. The first part of the proposition is thus proved. Now for each
$\lambda < 0$ close enough to $0$, there should exist some $t \in ]t_2 - \eta , t_2 + \eta [$
such that $L_t$ is tangent to $C^\lambda$. From what preceeds, $t$ has to be greater than $t_2$
and the proposition is proved, since $m(C^\lambda) = m(C) + 1$ when $\lambda < 0$. $\square$

\subsection{Neighbourhood of reducible curves}
\label{subsectreduc}

Let $C \in \R \overline{\cal M}_\gamma$ be a reducible curve and $C_1$, $C_2$ be its irreducible
components. For $i \in \{ 1,2 \}$, denote by $d_i = [C_i] \in H_2 (X ; \Z)$, $\underline{x}_i =
\underline{x} \cap C_i$ and $m_i = \# \underline{x}_i$. From Proposition \ref{propcompact},
$m_1 \in \{ c_1 (X) d_1 - 2 , c_1 (X) d_1 - 1 \}$. Denote by $t_3 = \pi_\gamma (C)$ and assume that
$\R C_1 \cap \R C_2 \cap L = \{ y \}$ and that $m_1 = c_1 (X) d_1 - 1$. Then, there exists
$\eta > 0$ such that the curves $C$ deforms to a one parameter family of real reducible
$J^t$-holomorphic curves $C^t_{red}$, $t \in ]t_3 - \eta , t_3 + \eta [$, which pass
through $\underline{x}^t$, where $(J^t , \underline{x}^t) = \gamma (t)$. The nodal point $y$
then deforms to a one parameter family of real non isolated double point $y^t$ of $\R C^t_{red}$.
Without loss of generality, we can assume that $y^t \notin B$ if $t \in ]t_3 - \eta , t_3 [$
and $y^t \in B$ if $t \in ]t_3 , t_3 + \eta [$.

\begin{prop}
\label{propred1}
Let $C^{t_3} = C_1^{t_3} \cup C_2^{t_3}  \in \R \overline{\cal M}_\gamma$ 
be a real reducible curve
and $t_3 = \pi_\gamma (C^{t_3})$. Assume that $\R C_1^{t_3} \cap \R C_2^{t_3} \cap L = \{ y^{t_3} \}$ 
and that 
$m_1 = c_1 (X) d_1 - 1$ with the above notations. Denote by $C^t_{red}$ (resp. $y^t$), 
$t \in ]t_3 - \eta , t_3 + \eta [$, the associated one parameter family of real reducible
$J^t$-holomorphic curves (resp. of real double point of $C^t_{red}$). Assume that $y^t \notin B$ 
if $t \in ]t_3 - \eta , t_3 [$ and $y^t \in B$ if $t \in ]t_3 , t_3 + \eta [$. Then,
as soon as $\eta$ is small enough, there exists a neighbourhood 
$W$ of $C$ in $\R \overline{\cal M}_\gamma$ such that for every $t \in ]t_3 - \eta , t_3 [$
(resp. $t \in ]t_3 , t_3 + \eta [$), $\sum_{C \in (\pi_\gamma^{-1} (t) \cap W)} <C,B> = -1$
(resp. $\sum_{C \in (\pi_\gamma^{-1} (t) \cap W)} <C,B> = +1$).
\end{prop} 
Note that all the curves $C^t$ close to $C^{t_3}$ are obtained topologically by smoothing the non 
isolated
real double point $y^{t_3}$ of $C^{t_3}$. Thus, they have the same mass as $C^{t_3}$.
$$\vcenter{\hbox{\input{rel4.pstex_t}}}$$
{\bf Proof:}

Without loss of generality, we can assume that $J^{t}$, $\underline{x}^{t}$ are constant and that $L$
(and the metric $g$) moves along a one parameter family $L_t$ which crosses the double point $y^{t_3}$ 
of $C^{t_3}$.
This indeed can be realized equivalently by fixing $L$ and having $J$, $\underline{x}$ moving
 along one parameter families $\phi_t^* J$, $\phi_t (\underline{x})$ where $\phi_t$ is some
$\Z/2\Z$-equivariant hamiltonian flow of $X$. The family of curves $C^t_{red}$ is then nothing but
the constant family $C$. Moreover, from Proposition $2.14$ of \cite{Wels1}, the curve $C^{t_3}$ 
extends to
a one parameter family $C^{t_3}_\lambda$, $\lambda \in ]-\epsilon , \epsilon [$, of real rational 
$J$-holomorphic curves which pass through $\underline{x}$ and realize $d$. These curves are obtained
topologically by smoothing the real double point $y^{t_3}$ of $C^{t_3}$. The intersection points 
between two
different curves in this family $(C^{t_3}_\lambda)_{\lambda \in ]-\epsilon , \epsilon [}$ are located
near the double points of $C^{t_3}$ and at $\underline{x}$. Thus, a neighbourhood $U$ of $y^{t_3}$ 
in $\R X$
is foliated by curves $C^{t_3}_\lambda \cap U$ and this foliation looks like the level
sets of an index one critical point of some Morse function $f : U \to \R$.
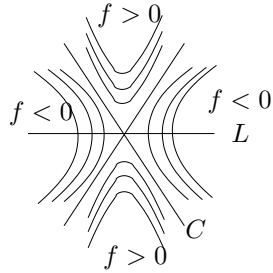
$$\vcenter{\hbox{\input{rel5.pstex_t}}}$$
We can assume that $L \cap U$ belongs to the domain $f \leq 0$. Let $(t_- , t_+) \in ]t_3 - \eta , t_3 [
\times ]t_3 , t_3 + \eta [$, restricting $U$ and $\epsilon$ if necessary, we can assume that
$L_{t_-}$ and $L_{t_+}$ are transversal to all the level sets $f \leq 0$. The number of maxima
minus the number of minima of $f$ restricted to $L_{t_\pm}$ is then equal to one, provided the latter
have been chosen generic. Now each maximum (resp. minimum) of $f$ restricted to $L_{t_-}$
corresponds to a curve $C^{t_-}$ having contact index $<C^{t_-} , B> = +1$ (resp. $<C^{t_-} , B> = -1$).
Likewise, each maximum (resp. minimum) of $f$ restricted to $L_{t_+}$
corresponds to a curve $C^{t_+}$ having contact index $<C^{t_+} , B> = -1$ (resp. $<C^{t_+} , B> = +1$),
hence the result. $\square$

\begin{prop}
\label{propred2}
Let $C^{t_4} = C_1^{t_4} \cup C_2^{t_4}  \in \R \overline{\cal M}_\gamma$ be a reducible curve
and $t_4 = \pi_\gamma (C^{t_4})$. Assume that $\R C_1^{t_4} \cap \R C_2^{t_4} \cap L = \{ y^{t_4} \}$ 
and that $m_1 = c_1 (X) d_1 - 2$ with the notations of Proposition \ref{propred1}. Then, there exist 
$\eta > 0$ and a neighbourhood $W$ of $C^{t_4}$ in $\R \overline{\cal M}_\gamma$ such that for 
every $t \in ] t_4 - \eta , t_4 + \eta [ \setminus \{ t_4 \}$, 
 $\sum_{C \in (\pi_\gamma^{-1} (t) \cap W)} <C,B> = 0$.
\end{prop} 
Note that once more, all the curves in $W$ have the same mass. Note also that $C_1^{t_4}$ belongs
to a one parameter family $C_1^{t_4} (\lambda)$ of $J^{t_4}$-holomorphic curves which pass through
$\underline{x}^{t_4}_1 = \underline{x}^{t_4} \cap C_1^{t_4}$ and realize $d_1$, whereas
$C_2^{t_4}$ does not deform to any $J^t$-holomorphic curve for $t \neq t_4$.\\

{\bf Proof:}

Without loss of generality, we can assume that $\underline{x}^{t}$ is constant. Let $U$ be a small
neighbourhood of $y^{t_4}$, it is foliated by the curves $C_1^{t_4} (\lambda) \cap U$. Choose a 
transversal $T$ to this foliation which is disjoint from $C_2^{t_4} \cap U$. From Proposition
$2.14$ of \cite{Wels1}, as soon as $\eta$ is small enough, there is one and only one
$J^t$-holomorphic real rational curve which pass through $\underline{x}^{t}$ and realize $d$
through every point of $T$. This produces a one parameter family of disjoint $J^t$-holomorphic real 
rational curve $\R C^t (\lambda) \cap U$, $\lambda \in T$. 
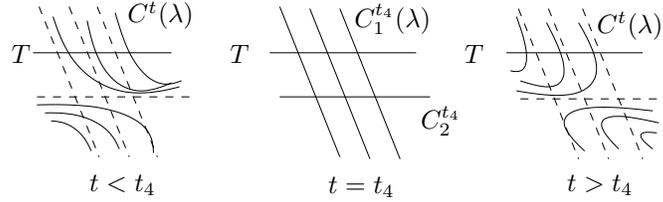
$$\vcenter{\hbox{\input{rel6.pstex_t}}}$$
Each of these curves $\R C^t (\lambda) \cap U$ has two connected components, which produce two 
functions partially defined on $L$ to $T$. To get the result, it is enough to observe that the 
number of maxima minus
the number of minima of these functions are either $+1$ and $-1$, or $0$ and $0$. $\square$

\begin{prop}
\label{propred3}
Let $C^{t_5} = C_1^{t_5} \cup C_2^{t_5} \in \R \overline{\cal M}_\gamma$ be a real reducible curve 
tangent to $L$ and $t_5 = \pi_\gamma (C^{t_5})$. Let $R$ be the number of real intersection points 
between $\R C_1^{t_5}$ and $\R C_2^{t_5}$. Then, there exist 
$\eta > 0$ and a neighbourhood $W$ of $C^{t_5}$ in $\R \overline{\cal M}_\gamma$ such that for 
every $t \in ] t_5 - \eta , t_5 + \eta [ \setminus \{ t_5 \}$, 
 $\pi_\gamma^{-1} (t) \cap W$ consists of exactly $R$ curves each of them obtained by smoothing
a different real intersection point between $\R C_1^{t_5}$ and $\R C_2^{t_5}$.
\end{prop} 

{\bf Proof:}

The proof is the same as the one of Proposition
$2.14$ of \cite{Wels1}, it is not reproduced here. The only argument which slightly differs from the
one in \cite{Wels1} is to show that for every real intersection point between $\R C_1^{t_5}$ and
$\R C_2^{t_5}$, there is at most one $J^t$-holomorphic curve in $\pi_\gamma^{-1} (t) \cap W$
obtained by smoothing this point. Actually, if there were two of them, they would intersect at
$\underline{x}^{t}$, at two points near each double point of $C^{t_5}$ but the one smoothed and
near the tangency point with $L$. This would produce more than $d^2$ intersection points, which
is impossible. $\square$

\subsection{Neighbourhood of the case $\protect\underline{x} \cap L \neq \emptyset$}
\label{subsectxL}

\begin{prop}
\label{propxL}
Let $C^{t_6} \in \R \overline{\cal M}_\gamma$ be an irreducible curve 
tangent to $L$ and $t_6 = \pi_\gamma (C^{t_6})$. Assume that $\underline{x}^{t_6}  \cap L = \{
x^{t_6}_1 \}$. Assume that $J^t$ is constant and that only the point $x^t_1$ actually depends on $t$.
Then, there exist 
$\eta > 0$ and a neighbourhood $W$ of $C^{t_6}$ in $\R \overline{\cal M}_\gamma$ such that for 
every $t \in ] t_6 - \eta , t_6 + \eta [ \setminus \{ t_6 \}$, 
 $\pi_\gamma^{-1} (t) \cap W$ consists of two curves having same mass and same contact index with
$L$ if $x^t_1$ is locally on the same size of $L$ as $C^{t_6}$, and $\pi_\gamma^{-1} (t) \cap W$
is empty otherwise.
\end{prop} 
$$\vcenter{\hbox{\input{rel7.pstex_t}}}$$

{\bf Proof:}

The moduli space of real rational $J^{t_6}$-holomorphic curves which pass through $\underline{x}^{t_6}
\setminus \{ x^{t_6}_1 \}$ and realize $d$ is one dimensional, and $C^{t_6}$ is a regular point
in this space. Thus, all the elements in this moduli space close to $C^{t_6}$ are located on
the same size of $L$ as $C^{t_6}$ itself. If $x^t_1$ is not on this size, we deduce that
$\pi_\gamma^{-1} (t) \cap W = \emptyset$ as soon as $W$ is small enough. Denote by $C^{t_6} (\lambda)$
the curves in this moduli space and let $U$ be a small neighbourhood of $x^{t_6}_1$ in $\R X$.
Then, the curves $(\R C^{t_6} (\lambda) \cap U) \setminus L$ have two connected components, which
produce two different foliations of one size of $L$ in $U \setminus L$ if $U$ is small enough. Thus, 
if $x^t_1$ is on this size, then $\# (\pi_\gamma^{-1} (t) \cap W) = 2$. In this case, the two curves
in $\pi_\gamma^{-1} (t) \cap W$ have obviously same mass and same contact index with $L$. $\square$

\subsection{Proofs of Theorems \ref{maintheo} and \ref{theorelation}}
\label{subsectionproofs} 

\subsubsection{Proof of Theorem \ref{maintheo}}
\label{subsubsectionproofmain}

Let $(J^0 , \underline{x}^0)$ and $(J^1 , \underline{x}^1)$ be two generic elements of
$\R {\cal J}_\omega \times \R_\tau X^{c_1 (X)d - 2}$ so that the integers $\Gamma_r^{d,B}
(J^0 , \underline{x}^0)$ and $\Gamma_r^{d,B} (J^1 , \underline{x}^1)$ are well defined.
Let $\gamma : t \in [0,1] \mapsto (J^t , \underline{x}^t) \in \R {\cal J}_\omega \times 
\R_\tau X^{c_1 (X)d - 2}$ be a generic path chosen in \S \ref{subsectchoice} joining
$(J^0 , \underline{x}^0)$ to $(J^1 , \underline{x}^1)$. Then, from genericity arguments
of \S \ref{subsectchoice}, we know that the integer $\Gamma_r^{d,B} (J^t , \underline{x}^t)$ 
is well defined for every
$t \in [0,1]$ but a finite number of parameters $0 < t_0 < t_1 < \dots < t_k < 1$ corresponding
to the following phenomena.

Concerning the first term in the definition of $\Gamma_r^{d,B} (J^t , \underline{x}^t)$:

1) Appearance of a unique real ordinary triple point or a unique real ordinary tacnode on an 
irreducible curve tangent to $L$.

2) Appearance of a transversal double point of an 
irreducible curve tangent to $L$ on $\underline{x}^t \cup L$.

3) Appearance of a real ordinary cusp of an 
irreducible curve on $L$.

4) Appearance of a an 
irreducible curve tangent to $L$ which is a critical point of $\pi_\gamma$ given by Theorem
\ref{theocrit}.

5) A sequence of curves of $\R {\cal M}_\gamma$ degenerates on a reducible curve given by 
Propositions \ref{propred1}, \ref{propred2} or \ref{propred3}.

6) One has $\underline{x}^t \cap L \neq \emptyset$.\\

Concerning the last three terms in the definition of $\Gamma_r^{d,B} (J^t , \underline{x}^t)$:

a) One of those considered in \cite{Wels2}.

b) A cuspidal curve has its cusp on $L$ but with a tangent line distinct from the one of $L$.

c) A reducible curve has one of the intersection points between its irreducible components on $L$
but is not tangent to $L$.

d) One has $\underline{x}^t \cap L \neq \emptyset$.\\

We have to prove that the integer $\Gamma_r^{d,B} (J^t , \underline{x}^t)$ does not change while
crossing one of these parameters $0 < t_0 < t_1 < \dots < t_k < 1$. In the cases 1, 2, a, this is
proven in the same way as in \cite{Wels1}, \cite{Wels2}. In the cases 3, b,
it follows from Proposition \ref{propcusp2}. Note that here the first term in the definition
of $\Gamma_r^{d,B} (J^t , \underline{x}^t)$ is not invariant. The term on cuspidal curves allows
to compensate for this lack of invariance. In the case 4, it follows from Propositions 
\ref{propinflectif}, \ref{propcusp1}. In the cases 5, c, it follows from Propositions \ref{propred1},
\ref{propred2} and \ref{propred3}. Note that here once more, in the case c, the first term
in the definition
of $\Gamma_r^{d,B} (J^t , \underline{x}^t)$ is not invariant. The term on reducible curves allows
to compensate for this lack of invariance. In the cases 6, d, it follows from Proposition 
\ref{propxL}. Here the first term in the definition
of $\Gamma_r^{d,B} (J^t , \underline{x}^t)$ is not invariant, this lack of invariance is compensated
thanks to the term on ${\cal T}an^d (J , \underline{x})$. $\square$

\subsubsection{Proof of Theorem \ref{theorelation}}
\label{subsubsectionproof}

Denote by $B (y , \epsilon)$ a disk of $\R X$ centered at $y \in X$ and having 
radius $\epsilon >0$.
Fix a generic $(J , \underline{x}) \in \R {\cal J}_\omega \times \R_\tau X^{c_1 (X)d - 2}$. When
$\epsilon$ converges to zero, $B (y , \epsilon) \to y$ and the three last terms in the definition
of $\Gamma_r^{d,B} (J , \underline{x})$ converge to $-\Gamma_r^{d} (J , \underline{x})$ since all the
curves does not move and all the special points are outside $B (y , \epsilon)$. At the same time,
the first term converges to a sum over real rational $J$-holomorphic curves which pass through
$\underline{x} \cup \{ y \}$ and realize $d$. Each of these curves are irreducible and immersed 
and deforms
in exactly two curves tangent to $\partial B (y , \epsilon)$ for $\epsilon \ll 1$. Moreover, the latter
are tangent from the outside of $B (y , \epsilon)$ and we deduce the relation $\Gamma_r^{d,B} =
2\chi^d_r - \Gamma_r^d$. 

Likewise, when $X = \C P^2$ and $\epsilon$ converges to $+\infty$, the three last terms 
in the definition
of $\Gamma_r^{d,B} (J , \underline{x})$ converge to $\Gamma_r^{d} (J , \underline{x})$. At the same time,
the first term converges to a sum over real rational $J$-holomorphic curves which pass through
$\underline{x}$, realize $d$ and are tangent to the line at infinity. Each of these curves are 
irreducible and immersed and deforms
in exactly two curves tangent to $\partial B (y , \epsilon)$ for $\epsilon \gg 1$. Moreover, 
one of these two curves is tangent from the outside of $B (y , \epsilon)$ and one from the inside, 
so that we get the relation $\Gamma_r^{d,B} = \Gamma_r^d$.

Finally, when $X = \C P^1 \times \C P^1$ and $\epsilon$ converges to $+\infty$, the boundary
of $B (y , \epsilon)$ accumulates on the union of a section $B_\infty$ and a fibre $F_\infty$ of
$\R P^1 \times \R P^1$. Then, the three last terms in the definition
of $\Gamma_r^{d,B} (J , \underline{x})$ converge to $\Gamma_r^{d} (J , \underline{x})$ as 
before. At the same time,
the first term converges to a sum over real rational $J$-holomorphic curves which pass through
$\underline{x}$, realize $d$ and are either tangent to $B_\infty \cup F_\infty$, or pass
through $B_\infty \cap F_\infty$. Each of the curves tangent to $B_\infty \cup F_\infty$ are 
irreducible and immersed and deforms
in exactly two curves tangent to $\partial B (y , \epsilon)$ for $\epsilon \gg 1$, one from the outside,
the other one from the inside. Likewise, each of the curves passing
through $B_\infty \cap F_\infty$ are irreducible and immersed and deforms
in exactly two curves tangent to $\partial B (y , \epsilon)$ for $\epsilon \gg 1$, both from the outside.
We hence get the relation $\Gamma_r^{d,B} =
2\chi^d_r + \Gamma_r^d$. $\square$

\section{On real conics tangent to five generic real plane conics}
\label{sectproofscon}

\subsection{Proofs of Theorem \ref{theocon} and Proposition \ref{propcomput}}
\label{subsectproofscon}

{\bf Proof of Theorem \ref{theocon}:}

The proof is similar to the one of Theorem \ref{maintheo}. We construct the universal moduli space
$\R {\cal C}_L$ of real pseudo-holomorphic conics tangent to $L_1, \dots , L_5$. It is a separable
Banach manifold of class $C^{l-k}$ equipped with a Fredholm projection $\pi_\R : \R {\cal C}_L \to
\R {\cal J}_\omega$ having vanishing index. Let $J_0$ and $J_1$ be two generic elements of
$\R {\cal J}_\omega$ so that $\Gamma^B (J_0)$ and $\Gamma^B (J_1)$ are well defined and
$\gamma : t \in [0,1] \mapsto J^t  \in \R {\cal J}_\omega $ be a generic path
joining $J_0$ to $J_1$. Denote by $\R {\cal C}_\gamma
= \R {\cal C}_L \times_\gamma [0,1]$, $\R \overline{\cal C}_\gamma$ its Gromov compactification
and $\pi_\gamma : \R \overline{\cal C}_\gamma \to [0,1]$ the associated projection. Genericity 
arguments similar to the ones of \S \ref{subsectchoice} show that the elements of 
$\R \overline{\cal C}_\gamma$ are smooth real conics
having a unique point of contact of order two with each $L_i$, $i \in \{ 1, \dots , 5 \}$, but
a finite number of them which may be:

1) Smooth real conics which are bitangent to $L_1 \cup \dots \cup L_5$, every
point of contact being of order at most two.

2) Smooth real conics which have a point of contact of order three with one curve $L_i$, 
$i \in \{ 1, \dots , 5 \}$, the other ones being non-degenerated.

3) Reducible conics made of two real lines, one of them being tangent to three curves $L_i$, 
$i \in \{ 1, \dots , 5 \}$, and the other one to the two remaining ones. These points of contact
are non-degenerated and outside the singular point of the conic.

4) Reducible conics of ${\cal C}on_{red}$ tangent to four curves $L_{i_1}, \dots , L_{i_4}$ and
whose singularity lie on the fifth curve $L_{i_5}$.

Likewise, the universal moduli space
$\R {\cal C}_L^{red}$ of real reducible pseudo-holomorphic conics tangent to four curves
$L_{i_1}, \dots , L_{i_4}$ out of the five $L_1, \dots , L_5$ is a separable
Banach manifold of class $C^{l-k}$. Denote by $\R {\cal C}_\gamma^{red}
= \R {\cal C}_L^{red} \times_\gamma [0,1]$. It is a one dimensional compact manifold whose
elements are couples of real lines having four points of contacts with $L_{i_1}, \dots , L_{i_4}$ 
which are of order two, but a finite number of them which may be:

a) Tangent to the five curves $L_1, \dots , L_5$, with non-degenerated points of contacts.

b) Tangent to $L_{i_1}, \dots , L_{i_4}$ but with one point of contact of order three.

c) Tangent to $L_{i_1}, \dots , L_{i_4}$ with their singular point on $L_{i_5}$.

d) Tangent to $L_{i_1}, \dots , L_{i_4}$ with their singular point on $L_{i_1} \cup \dots \cup L_{i_4}$.

The only thing to check is that the value of $\Gamma^B (J_t)$ does not change while $t$ crosses one
of the special values listed in 1-4 and a-d. In the cases $1$, a and d, it is easy to check. In the
cases $2$, $b$, it follows from Proposition \ref{propinflectif}. In the case $3$, the proof is the same
as the one of Proposition \ref{propred3}. Finally, in the cases $4$, c, the proof is the same
as the one of Proposition \ref{propred1}. $\square$\\

{\bf Proof of Proposition \ref{propcomput}:}

From Theorem \ref{theocon}, we can assume that the five disjoint disks are of radius $\epsilon$
small, and have $\epsilon$ converging to zero so that they contracts onto five distinct
points $y_1, \dots, y_5$. The conics tangent to $L_1, \dots , L_5$ degenerate onto conics passing
through $y_1, \dots, y_5$. From \cite{Gro}, there is only one such $J$-holomorphic conic. Reversing
this process as in the proof of Theorem \ref{theorelation}, each conic passing through $y_i$ 
deforms into two conics which are tangent to $B_i$ from the
outside, for $\epsilon$ small enough. As soon as $\epsilon$ is small enough, the first term in 
the definition of $\Gamma^B$ equals then $2^5=32$. Likewise, elements of ${\cal C}on_{red}$
degenerate onto reducible conics passing
through four out of the five points $y_1, \dots, y_5$. There are five ways to choose these four points,
three couples of lines passing through these four points and each of these couples deforms into
$2^4=16$ reducible conics tangent to the four associated disks $B_i$ from the outside, as soon as 
$\epsilon$ is small enough. Since the singular point of these conics is outside $B = B_1 \cup \dots
\cup B_5$, the second term in 
the definition of $\Gamma^B$ equals $-5*3*16 = -240$. We deduce that $\Gamma^B = 32+240 = 272$.
Likewise, if $B_1, \dots , B_5$ are close to five generic double lines of the plane we can have
the curves $L_1, \dots , L_5$ degenerate onto five couples of real lines $L_i^1 \cup L_i^2$ close
to the double lines and intersecting each other at $x_1, \dots, x_5$. Every conic tangent to
$L_i$ degenerates onto a conic tangent to $L_i^1 \cup L_i^2$ or a conic which passes through
$x_i$. Now each conic tangent to $L_i^1$ deforms to a conic tangent to $L_i^2$ since $L_i^1$ and
$L_i^2$ are as close to each other as we wish. Hence these conics come by pairs, one deforming
to a conic tangent from the inside of $B_i$ and the other one from the outside. Hence, the only conics
which contributes to the first term of $\Gamma^B$ correspond to the ones passing through
$x_1, \dots, x_5$. Their contribution is $32$ as before. In the same way, the second term
 of $\Gamma^B$ equals $-240$ as before, as soon as $B_1, \dots , B_5$ are close enough to
$(L_1^1 \cup L_1^2), \dots, (L_5^1 \cup L_5^2)$. Hence the result. $\square$

\subsection{How does $\Gamma^B$ depend on the isotopy class of $B$?}
\label{subsectdepend}

Let $B_2, \dots , B_5$ be four disks of $\R P^2$ transversal to each other and $(B_1^t)_{t \in
]-\epsilon, \epsilon[}$ be a smooth one-parameter family of disks which are transversal to
$B_2, \dots , B_5$ for $t \in ]-\epsilon, \epsilon[ \setminus \{ 0 \}$ and which have an order two
point of contact $x$ with $B_2$ for $t=0$. We can assume that for $t \in ]-\epsilon, 0[$ (resp.
for $t \in ] 0 , \epsilon[$), the curves $L_1^t = \partial B_1^t$ and $L_2 = \partial B_2$ have
two intersection points (resp. do not intersect) in a neighbourhood of $x$.
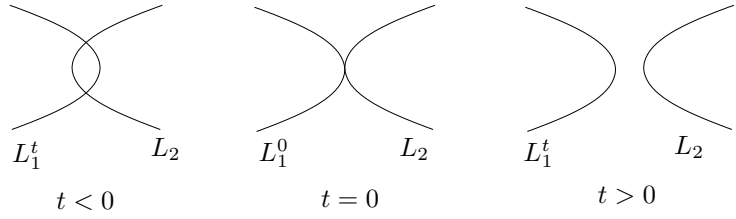
$$\vcenter{\hbox{\input{rel8.pstex_t}}}$$
Denote by $B^t = B_1^t \cup B_2 \cup \dots \cup B_5$, the integer $\Gamma^{B^t}$ is well defined
for $t \in ]-\epsilon, \epsilon[ \setminus \{ 0 \}$. We have to compare the values of $\Gamma^{B^t}$
for $t<0$ and $t>0$. Denote by ${\cal C}on (J , x)$ the finite set of real conics which are tangent
to $L_3 , L_4 , L_5$, pass through $x$ and are tangent at $x$ to $L_1^0$ and $L_2$. Likewise, 
denote by ${\cal C}on_{red} (J , x)$ the finite set of real reducible conics made of the $J$-holomorphic
line $T_x$ which is tangent to $L_1^0$ and $L_2$ at $x$ and of a real $J$-holomorphic line 
tangent to two curves out of the three curves $L_3 , L_4 , L_5$.
\begin{prop}
\label{propdepend}
Let $B^t = B_1^t \cup B_2 \cup \dots \cup B_5$ be a one parameter family of five disks in $\R P^2$
as above and $(t_- , t_+) \in ]-\epsilon, 0[ \times ] 0 , \epsilon[$. Then,
$$\Gamma^{B^{t_+}} = \Gamma^{B^{t_-}} + 2<L_2 ,B_1^0> \Big( \sum_{C \in {\cal C}on (J , x)} \Pi_{j=2}^5 
<C , B_j> \Big) -
2<L_2 ,B_1^0> \Big( \sum_{C \in {\cal C}on_{red} (J , x)} \Pi_{j=2}^5 <C , B_j> \Big).$$
\end{prop}
{\bf Proof:}

We can have locally $L_2$ degenerate on a half line. The conics tangent to $L_2$ degenerate then on 
conics tangent to the half line and conics passing through the vertice $s$ of this half line. As in the 
proof of Theorem \ref{theorelation}, the contribution to $\Gamma^{B^t}$ of conics tangent to the 
half line vanishes. As $t$ goes to zero, $s$ converges to $x$ and the conics passing through $s$
and tangent to $B_1^t, B_3, B_4, B_5$ converge to conics passing through the vertice $s$ and tangent
to $B_1^0, B_3, B_4, B_5$. If the order two point of contact of the latter with $B_1^0$ is outside
$x$, they can be deformed for $t \in ]-\epsilon, \epsilon[$. If on the contrary such a conic $C$ belongs
to ${\cal C}on (J , x) \cup {\cal C}on_{red} (J , x)$, it follows from Proposition \ref{propxL}
that it deforms for $t \in ]-\epsilon, 0]$ (resp. $t \in [ 0, \epsilon[$) if and only if
$<C,B_1^0> = -<L_2 ,B_1^0>$ (resp. $<C,B_1^0> = <L_2 ,B_1^0>$), that is, if and only if $C$ and $L_2$
are locally on opposite sides of $B_1^0$ (resp. on the same side of $B_1^0$). Hence the result.
$\square$

\subsection{Final remarks}
\label{subsectfinal}

1) The results of \S \ref{sectproofscon} take advantage of the fact that a pseudo-holomorphic conic
cannot be cuspidal and may have two irreducible components at most. To extend the results of
\S \ref{subsectresults} to pseudo-holomorphic curves having $s > 1$ tangency conditions with $L$
would seem to require the introduction of $4^s$ terms in the definition of $\Gamma^{d,B}$. These terms
consist of curves having $s_1$ tangency conditions with $L$, $s_2$ cusps, $s_3 + 1$ irreducible 
components and $s_4$ tangency conditions with the lines $T_i$, $i \in I$, where $s_1 + \dots + s_4 = s$.
One should then study the collisions between these tangency conditions, cusps, etc... which has not
been done here.

2) In contrast with the works \cite{Wels1}, \cite{Wels2}, the moduli space $\R {\cal M}^d_L$ do not
appear here as the fixed point set of some $\Z/2\Z$-action on some complexified moduli space
${\cal M}^d_L$. For such a purpose, we should have complexified $L$ to some surface $L_\C$
in $X$ and restricted ourselves to almost complex structures $J$ for which $L_\C$ is $J$-antiholomorphic,
as in \cite{LiRu} and \cite{IoPar}. The advantage not to do so here was to get immediatly some
invariant for any $J \in \R {\cal J}_\omega$ without any restrictions.

3) The condition that $L$ is smooth and bounded by a smooth surface $B$ is of course too restrictive.
We reduced our study to this case for convenience. For example, one could replace the embedding 
$B \to \R X$ with some smooth map with finitely many ramification points and which maps the boundary $L$
of $B$ to some immersed curve with transversal double points as singularities. The index $<x , B>$
for $x \in \R X$ should then be defined as twice the number of preimages of $x$ in $B$ less one.
Since every step of the proof of Theorem \ref{maintheo} is local, it readily extends to this
case.

\addcontentsline{toc}{part}{\hspace*{\indentation}Bibliography}


\bibliography{relative}
\bibliographystyle{abbrv}

\noindent \'Ecole normale sup\'erieure de Lyon\\
Unit\'e de math\'ematiques pures et appliqu\'ees\\
UMR CNRS $5669$\\
$46$, all\'ee d'Italie\\
$69364$, Lyon cedex $07$\\
(FRANCE)\\
e-mail : {\tt jwelschi@umpa.ens-lyon.fr}

\end{document}